\documentclass[a4paper]{amsart}
\usepackage{amsmath}
\usepackage{amscd}
\usepackage{a4wide}               
\usepackage{comment,url}
\usepackage[pdfpagelabels,
colorlinks=true,citecolor=blue,hypertexnames=false]{hyperref}

\subjclass[2000]{Primary 05A15, 14N10, 33F10, 68W30; Secondary 33C05, 97N80}

\newtheorem{thm}{Theorem}
\newtheorem{cor}[thm]{Corollary}
\newtheorem{lem}[thm]{Lemma}

\def\<#1>{\langle#1\rangle}
\let\set\mathbb
\def\deg{\operatorname{deg}}
\def\ord{\operatorname{ord}}
\def\res{\operatorname{res}}

\def\cand{\mathrm{cand}}

\begin{document}

 \title{The complete~generating~function for~Gessel~walks is~algebraic}
\author[Alin Bostan]{Alin Bostan}
\address{Algorithms Project, INRIA Paris-Rocquencourt, 78153 Le  
Chesnay, France}
\email{Alin.Bostan@inria.fr}

\author[Manuel Kauers]{Manuel Kauers}
\address{Research Institute for Symbolic Computation, J. Kepler  
University Linz, Austria}
\email{mkauers@risc.uni-linz.ac.at}
\keywords{Combinatorial enumeration,  generating function, lattice walks, Gessel conjecture, algebraic functions, computer algebra, 
automated guessing, fast algorithms.}

\begin{abstract}
Gessel walks are lattice walks in the quarter plane $\set N^2$ which start at the
origin~$(0,0)\in\set N^2$ and consist only of steps chosen from the set
$\{\leftarrow,\swarrow,\nearrow,\rightarrow\}$. We prove that if $g(n;i,j)$ denotes the
number of Gessel walks of length~$n$ which end at the point~$(i,j)\in\set N^2$, then the
trivariate generating series $ \displaystyle{ G(t;x,y)=\sum_{n,i,j\geq 0} g(n;i,j)x^i y^j
t^n } $ is an algebraic function.
\end{abstract}

 \maketitle

 \section{Introduction}

 The starting question in lattice path theory is the following:
 How many ways are there to walk from the origin through the lattice $\set Z^2$ to a specified
 point~$(i,j)\in\set Z^2$, using a fixed number~$n$ of steps chosen from a given set~$S$
 of admissible steps. The question is not hard to answer. If we write $f(n;i,j)$
 for this number and define the generating function 
 \[
  F(t;x,y):=\sum_{n=0}^\infty\Bigl(\sum_{i,j\in\set Z}f(n;i,j)x^i y^j\Bigr)t^n
\; 
\in \set Q[x,y,x^{-1},y^{-1}][[t]]
 \]
 then a simple calculation suffices to see that $F(t;x,y)$ is rational, i.e., it agrees with the series expansion at $t=0$
 of a certain rational function $P/Q\in\set Q(t,x,y)$.
 This is elementary and well-known.

 Matters are getting more interesting if restrictions are imposed. For example,
 the generating function $F(t;x,y)$ will typically no longer be rational if lattice
 paths are considered which, as before, start at the origin, consist of $n$ steps,
 end at a given point~$(i,j)$, but which, as an additional requirement, 
 \emph{never step out of the right half-plane.} 
 In was shown in~\cite[Prop.~2]{melou03} 
that no matter which set $S$ of admissible steps
 is chosen, the complete generating function $F$ for such walks is algebraic, 
 i.e., it satisfies $P(F,t,x,y)=0$ for some polynomial $P\in\set Q[T,t,x,y]$.

 If the walks are not restricted to a half-plane but to a quarter plane, say to
 the first quadrant, then the generating function $F$ might not even be algebraic. 
 For some step sets it is, for others it is not~\cite{bousquet05,mishna09}. Among the series which are not 
 algebraic, there are some which are still D-finite with respect to~$t$ (i.e., they satisfy a linear 
 differential equation in~$t$ with polynomial coefficients in $\set Q[t,x,y]$), and others which are not
 even that~\cite{melou03,rechnitzer07}.

 Bousquet-M\'elou and Mishna~\cite{BoMi08} have systematically investigated all the walks in the
 quarter plane with step sets $S\subseteq\{{\leftarrow,}\ {\nwarrow,}\ {\uparrow,}\
 {\nearrow,}\ {\rightarrow,}\ {\searrow,}\ {\downarrow,}\ {\swarrow}\}$. 
 After discarding trivial cases and applying
 symmetries, they reduced the 256 different step sets to 79 inherently different
 cases to study. They provided a unified way to prove that 22 of those are D-finite,
 and gave striking evidence that 56 are not D-finite. Only a single step set 
 sustained their attacks, and this is the step set that we are considering here. 

 This critical step set is $\{\leftarrow,\swarrow,\nearrow,\rightarrow\}$. 
 The central object of the present article are thus lattice walks in $\set Z^2$ which
 \begin{itemize}
 \item start at the origin $(0,0)$,
 \item consist of $n$  steps chosen from the step set
   $\{\leftarrow,\swarrow,\nearrow,\rightarrow\}$, and
 \item never step out of the first quadrant $\set N^2$ of $\set Z^2$. 
 \end{itemize}
 These walks are also known as \emph{Gessel walks}.
 By $g(n;i,j)$, we denote the number of Gessel walks of length~$n$ which
 end at the point~$(i,j)\in\set Z^2$. The complete generating function 
 of this sequence is denoted by
 \[
  G(t;x,y)=\sum_{n=0}^\infty\Bigl(\sum_{i,j\in\set Z}
    g(n;i,j)x^i y^j\Bigr)t^n.
 \]
 Since $g(n;i,j)=0$ if $\min(i,j)>n$ or $\max(i,j)<0$, the inner sum is a polynomial 
 in $x$ and~$y$ for every fixed choice of~$n$, and thus $G(t;x,y)$ lives in $\set Q[x,y][[t]]$.

 Gessel [unpublished] considered the special end point $i=j=0$, i.e., 
 Gessel walks returning to the origin, so-called \emph{excursions}. Their counting sequence $g(n;0,0)$
 starts as
 \[
  1,\ 0,\ 2,\ 0,\ 11,\ 0,\ 85,\ 0,\ 782,\ 0,\ 8004,\ 0,\ 88044,\ 0,\ 1020162,\ 0,\ \dots
 \]
 He observed empirically that these numbers admit a simple hypergeometric 
 closed form. His observation became known as the \emph{Gessel
 conjecture}, and remained open for several years. Only recently, it was shown to
 be true:

\begin{thm}\label{thm:0} \cite{kauers08g} $\displaystyle G(t;0,0)=
{_3F_2}\biggl(\begin{matrix}
  {5/6}\kern.707em {1/2}\kern.707em{1}\\ 
  {5/3}\kern1em{2}
\end{matrix}\,\bigg|\,16t^2\biggr)
=\sum_{n=0}^\infty\frac{(5/6)_n(1/2)_n}{(5/3)_n(2)_n}(4t)^{2n}$.
\end{thm}

This result obviously implies that $G(t;0,0)$ is D-finite. 
Less obvious at this point, and actually overlooked until now, is the
fact that the power series $G(t;0,0)$ is even \emph{algebraic.} 
Because of the alternative representation
\[
{_3F_2}\biggl(\begin{matrix}
  {5/6}\kern.707em {1/2}\kern.707em{1}\\ 
  {5/3}\kern1em{2}
\end{matrix}\,\bigg|\,16t^2\biggr)
= \frac{1}{t^2} \left(\frac12\, {_2F_1}\biggl(\begin{matrix}
  {-1/6}\kern.707em {-1/2}\\ 
  {2/3}\kern1em{}
\end{matrix}\,\bigg|\,16t^2\right) - \frac12\biggr)
\]
it was clear that algebraicity could be decided by reference to
Schwarz's classification~\cite{schwarz1872} of algebraic~$_2F_1$'s, 
but as nobody recognized that the parameters $(-1/6,-1/2;2/3)$ 
actually fit to Case~III of Schwarz's table, the rumor started to 
circulate that $G(t;0,0)$ is not algebraic. In fact:  

\begin{cor}\label{cor:2} $G(t;0,0)$ is algebraic.
\end{cor}                             
With Theorem~\ref{thm:0} and
standard software packages like \textsf{gfun}~\cite{salvy94,mallinger96} at hand, discovering and
proving Cor.~\ref{cor:2} is an easy computer algebra exercise.
Compared to a proof by table-lookup, the constructive proof given
below has the advantage that it applies similarly also for families
of functions for which classification results are not available.

\begin{proof}    
The idea is to come up with a polynomial $P(T,t)$ in $\mathbb{Q}[T,t]$ 
and prove that $P$ admits the power series $g(t) =
 \sum_{n=0}^\infty\frac{(5/6)_n(1/2)_n}{(5/3)_n(2)_n}(16t)^{n}$ as a root.
Using Thm.~\ref{thm:0}, this implies that $P(T,t^2)$ is an annihilating polynomial for $G(t;0,0)$, so that the latter series is indeed algebraic.

Such a polynomial $P$ can be \emph{guessed\/} starting from the first, say, 100 terms, of the series $g(t)$, using for instance \textsf{Maple}'s
routine \textsf{seriestoalgeq} from the \textsf{gfun} package (see Sections~\ref{sec:2.1} and~\ref{sec:3.1} for more details on automated guessing). The explicit form of $P$ is given below.
               
By the implicit function theorem, that polynomial $P$ admits a root $r(t) \in \set Q[[t]]$ with $r(0)=1$.
Since $P(T,0)=T-1$ has a single root in $\set C$, the series $r(t)$ is the unique root of $P$ in~$\set C[[t]]$. 
Now, $r(t)$ being algebraic, it is D-finite, and thus its coefficients satisfy a recurrence with polynomial coefficients. To complete the proof, it is then sufficient to type the following commands into \textsf{Maple}.
 \begin{verbatim}
> with(gfun):
> P:=(t,T) -> -1+48*t-576*t^2-256*t^3+(1-60*t+912*t^2-512*t^3)*T+(10*t
   -312*t^2+624*t^3-512*t^4)*T^2+(45*t^2-504*t^3-576*t^4)*T^3+(117*t^3
   -252*t^4-288*t^5)*T^4+189*t^4*T^5+189*t^5*T^6+108*t^6*T^7+27*t^7*T^8:
> gfun:-diffeqtorec(gfun:-algeqtodiffeq(P(t,r), r(t)), r(t), g(n));
\end{verbatim}
This outputs the first-order recurrence
\[
(n + 2) (3 n + 5) g_{n + 1} - 4 (6 n + 5) (2 n + 1) g_n = 0, \qquad g_0=1, 
\]          
satisfied by the coefficients of $r(t) = \sum_{n=0}^\infty g_n t^n$. 
Its solution is $g_n=\frac{(5/6)_n(1/2)_n}{(5/3)_n(2)_n}16^n$, and therefore $g(t)$ and $r(t)$ coincide, and thus $g(t)$ is a solution of $P$, as was to be shown.
\end{proof}

\bigskip
The aim in the present article is to lift the result of Corollary~\ref{cor:2}
to the complete
generating function, where $x$ and $y$ are kept as parameters. We are going to show:

\begin{thm}\label{thm:g} $G(t;x,y)$ is algebraic.
\end{thm}

This twofold generalization of Thm.~\ref{thm:0} is a surprising result.
Until now, it was not known whether $G(t;x,y)$ is even D-finite with respect to~$t$ or not, 
and both cases seemed equally plausible in view of known results about other step sets. Thm.~\ref{thm:g} implies that $G(t;x,y)$ is D-finite with respect to each of its variables, and in particular that the sequence $g(n;i,j)$ is P-finite
(i.e., it satisfies a linear recurrence with polynomial coefficients in~$n$)
for any choice of $(i,j) \in \set N^2$. This settles several conjectures made by Petkov{\v s}ek and Wilf in~\cite[\S2]{petkovsek08}.   
As noted in~\cite{petkovsek08}, even for simple values of $(i,j)$ the sequence $g(n;i,j)$ is not hypergeometric, unlike the excursions sequence $g(2n;0,0)$.   
For instance, the sequence $g(2n+1;1,0)$ satisfies a third order linear recurrence, 
but it is not hypergeometric. 
Moreover, no closed formula seems to exist for $g(n;i,j)$, for arbitrary $(i,j)$. 
All this indicates that counting general walks is much more difficult that just counting excursions.


\medskip Theorem~\ref{thm:g} will be established by obstinately using the approach based on \emph{automatic guessing and proof\/} promoted in~\cite{bostan09}, and by making heavy use of computer algebra.
In contrast to Corollary~\ref{cor:2}, we manage in our proof of Theorem~\ref{thm:g} to avoid exhibiting a polynomial that has $G(t;x,y)$ as a root.
This is fortunate, since a posteriori estimates show that the minimal polynomial of $G(t;x,y)$ is huge, having a total size of about 30Gb.

Only annihilating polynomials of the section series $G(t;x,0)$ and $G(t;0,y)$
are produced and manipulated during the computer-driven proof of Theorem~\ref{thm:g}.
But even restricted to those ones,
our computations have led to expressions far too large to be included into a 
printed publication; too large even to be processed efficiently by standard
computer algebra systems like \textsf{Maple} or \textsf{Mathematica}. To get the computations completed,
it was necessary to use careful implementations of sophisticated special 
purpose algorithms, and to run these on computers equipped with fast processors 
and large memory capacities. These computations were performed using the computer algebra system \textsf{Magma}~\cite{magma}.
Our result is therefore interesting not only because of its combinatorial significance,
but it is also noteworthy because of the immense computational effort that was 
deployed to establish it. 

\section{A Dry Run: Kreweras walks}\label{sec:kreweras}

The computations which were needed for proving Thm.~\ref{thm:g} were performed by means of
efficient special purpose software running on fast hardware. It would not be easy to redo
these calculations in, say, \textsf{Maple} or \textsf{Mathematica} on a standard 
computer.
As a more easily 
reproducible calculation, we will show in this section how to reprove the classical result 
that the generating function of Kreweras walks is algebraic~\cite{kreweras65,gessel86,bousquet05}.
A slight variation of the very same reasoning, albeit with intermediate expressions far too large to be spelled out
here, is then used in the next section to establish Thm.~\ref{thm:g}.

Kreweras walks differ from Gessel walks only in their choice of admissible steps.
They are thus defined as lattice walks in~$\set Z^2$ which
\begin{itemize}
\item start at the origin $(0,0)$,
\item consist only of steps chosen from the step set
  $\{\leftarrow,\downarrow,\nearrow\}$, and
\item never step out of the first quadrant $\set N^2$ of $\set Z^2$. 
\end{itemize}
If $f(n;i,j)$ denotes the number of Kreweras walks consisting of $n$ steps and ending 
at the point $(i,j)\in\set Z^2$, then it follows directly from 
its combinatorial definition
that the sequence $f(n;i,j)$ satisfies the multivariate recurrence 
with constant coefficients
 \begin{equation}\label{recf}
 f(n+1;i,j) = f(n;i+1,j) + f(n;i,j+1) + f(n;i-1,j-1),
\end{equation}
for all $n,i,j\geq 0$.
Together with the boundary conditions $f(n;-1,0)=f(n;0,-1)=0$ ($n\geq0$) and $f(0;i,j)=\delta_{i,j,0}$
($i,j\geq0$), this recurrence equation implies the functional equation
\[
   F(t;x,y)=1 + \bigl(\tfrac1x+\tfrac1y+xy\bigr)t F(t;x,y)
         - \tfrac1y t F(t;x,0) - \tfrac1x t F(t;0,y)
\]
for the generating function 
\[
 F(t;x,y)=\sum_{n=0}^\infty\Bigl(\sum_{i,j=0}^\infty f(n;i,j)x^iy^j\Bigr)t^n.
\]
Noting that $F(t;0,y)$ and $F(t;y,0)$ are equal by the symmetry of the step set about the main diagonal of $\set N^2$, the last equation becomes
\[ 
 F(t;x,y)=1 + \bigl(\tfrac1x+\tfrac1y+xy\bigr)t F(t;x,y)
         - \tfrac1y t F(t;x,0) - \tfrac1x t F(t;y,0).
\]

At the heart of our next arguments is the \emph{kernel method}, a method commonly attributed to Knuth~\cite[Solutions of Exercises 4 and 11 in \S2.2.1]{knuth69} 
which has already been used to great advantage in lattice path counting, see e.g.~\cite{fayolle79,prodinger03,bousquet05}. After bringing the
functional equation for $F(t;x,y)$ to the form
\begin{align}\label{eq:K}\tag{K} 
 ((x+y+x^2y^2)t-xy)F(t;x,y) &= x t F(t;x, 0) + y t F(t;y,0) - x y,
\end{align}
the kernel method consists of coupling $x$ and $y$ in such a way that this equation 
reduces to a simpler one, from which useful information about the \emph{section series\/} $F(t;x,0)$ can be 
extracted. In the present case, the substitution
\begin{alignat*}1 
 y \rightarrow Y(t,x)
   &=\frac{x-t-\sqrt{-4 t^2 x^3+x^2-2 t x+t^2}}{2 t x^2}\\
   &=t + \tfrac1x t^2 + \tfrac{x^3+1}{x^2}t^3 + \tfrac{3x^3+1}{x^3}t^4 + \tfrac{2x^6+6x^3+1}{x^4}t^5 + \cdots      \in \set Q[x,x^{-1}][[t]], 
\end{alignat*}
which is legitimate since the power series $Y(t,x)$ has positive valuation,
puts the left hand side of~\eqref{eq:K} to zero, and therefore shows that $U = F(t;x,0)$ is a solution of the \emph{reduced kernel equation\/} 
\begin{equation}\label{eq:Kred}\tag{$\mathrm{K}_\mathrm{red}$}
   U(t,x) = \frac{Y(t,x)}{t} - \frac{Y(t, x)}{x} U(t,Y(t,x)).
\end{equation}                         

\medskip
Now, the key feature of Equation~\eqref{eq:Kred} is that its unique solution in~$\set Q[[x,t]]$ is $U=F(t;x,0)$. This is a consequence of the following easy lemma. Here, and in the rest of the article, $\ord_v S$ denotes the valuation of a power series $S$ with respect to some variable $v$
occurring in $S$.
\begin{lem}\label{lem:1}
 Let $A,B,Y\in\set Q[x,x^{-1}][[t]]$ be such that  $\ord_t B>0$ and $\ord_t Y>0$. 
 Then there exists at most one power series $U\in\set Q[[x,t]]$ with
 \[
   U(t,x) = A(t, x) + B(t, x) \cdot U(t, Y(t, x)).
 \]
\end{lem}
\begin{proof}                                                 
By linearity, it suffices to show that the only solution in $\set Q[[x,t]]$ of the homogeneous equation $U(t,x) = B(t, x) \cdot U(t, Y(t, x))$ is the trivial solution $U=0$. This is a direct consequence of the fact that if $U$ were non-zero, then the valuation of $B(t, x) \cdot U(t, Y(t, x))$ would be 
at least equal to $\ord_t B + \ord_t U$,
thus strictly greater than the valuation of $U(t,x)$, a contradiction.
\end{proof}
          
We are now ready to reprove the following classical result.

\begin{thm}\label{thm:k}\cite{gessel86} $F(t;x,y)$ is algebraic. 
\end{thm}

\begin{proof} The strategy is to use a computer-assisted proof, 
which is completed in two steps:
\begin{enumerate}
\item\label{step:1} \emph{Guess} an algebraic equation for the series $F(t;x,0)$, by inspection  of its initial terms.
\item\label{step:2} \emph{Prove} that 
\begin{enumerate}   
	\item\label{step:2a}   
	the equation guessed at Step~\eqref{step:1} admits exactly one solution in $\set Q[[x,t]]$, denoted $F_{\cand}(t;x,0)$;
	\item\label{step:2b} 
	the power series $U = F_{\cand}(t;x,0)$ satisfies~\eqref{eq:Kred}.    
\end{enumerate}
\end{enumerate}
Once this has been accomplished, the fact that 
$U = F(t;x,0)$ also satisfies Equation~\eqref{eq:Kred}, in conjunction with Lemma~\ref{lem:1} (with the choice $A(t,x) = Y(t,x)/t$ and $B(t,x)=-Y(t,x)/x$), implies that the power series $F_{\cand}(t;x,0)$ and $F(t;x,0)$ coincide.

In particular, $F(t;x,0)$ satisfies the guessed equation, and 
this certifies that $F(t;x,0)$ is an algebraic power series. Since $Y(t,x)$ is algebraic as well, and since the class of algebraic power series is closed under addition, multiplication and inversion,
it follows from \eqref{eq:K} that $F(t;x,y)$ is algebraic, too. This  concludes the proof. 
\end{proof}

In the rest of this section, we supply full details on the automated guessing step~\eqref{step:1} and on the proving steps~\eqref{step:2a} and~\eqref{step:2b}. 
                                                  
\subsection{Guessing}\label{sec:2.1}

Given the first few terms of a power series, it is possible to determine potential equations that the power series may satisfy, for example by making a suitable ansatz with undetermined coefficients and solving a linear system. In practice, either Gaussian elimination, or faster, special purpose algorithms based on Hermite-Pad\'e approximation~\cite{beckermann94}, are used. The computation of such candidate 
equations is known as \emph{automated guessing\/} and is one of the most widely known features of 
packages such as \textsf{Maple}'s \textsf{gfun}~\cite{salvy94}. 

If sufficiently many terms of the series are provided, automated guessing will 
eventually find an equation whenever there is one. 
The method has two possible drawbacks. First, it may in principle return
false equations (although, if applied properly, it virtually never does so
in practice). This is why -- in order to provide fully rigorous proofs -- equations discovered by this method must be subsequently proven
by an independent argument.
Second, if the precision needed to recover the equations is very high, the guessing computations could take extremely long when using traditional software.
This is typically the case in the Gessel example treated in Section~\ref{sec:gessel}, for which dedicated, very efficient, algorithms are needed. 

\medskip In the Kreweras case, the computations are feasible in \textsf{Maple}.
We now provide commented code which enables the discovery of 
an algebraic equation potentially satisfied by~$F(t;x,0)$.
First, a function $f$ is defined which computes the numbers $f(n;i,j)$ 
via the multivariate
recurrence~\eqref{recf}.
\begin{verbatim}  
> f:=proc(n,i,j) 
  option remember; 
    if i<0 or j<0 or n<0 then 0 
    elif n=0 then if i=0 and j=0 then 1 else 0 fi 
    else f(n-1,i-1,j-1)+f(n-1,i,j+1)+f(n-1,i+1,j) fi 
  end:    
\end{verbatim} 

Using this function, we compute the first 80 coefficients of $F(t;x,0)$;
they are polynomials in $x$ with integer coefficients. The resulting truncated power series is stored in the variable $S$.    
\begin{verbatim}
> prec:=80:
> S:=series(add(add(f(k,i,0)*x^i,i=0..k)*t^k,k=0..prec),t,prec-1):
\end{verbatim}

Next, starting from $S$, the \textsf{gfun} guessing function \textsf{seriestoalgeq} discovers a candidate for an algebraic equation satisfied by $F(t;x,0)$. For efficiency reasons, 
we do not use the built-in version of \textsf{gfun}, but a recent one which can be
downloaded from \url{http://algo.inria.fr/libraries/papers/gfun.html}

\begin{verbatim}
> gfun:-seriestoalgeq(S,Fx(t)):
> P:=collect(numer(subs(Fx(t)=T,%[1])),T); 
\end{verbatim} 
The guessed polynomial reads:            
\begin{alignat*}1   \label{guessedP}
 P(T,t,x)&=
   (16 x^3 t^4+108 t^4-72 x t^3+8 x^2 t^2-2 t+x)
  \\&\quad+(96 x^2 t^5-48 x^3 t^4-144 t^4+104 x t^3-16 x^2 t^2+2 t-x)T
  \\&\quad+(48 x^4 t^6+192 x t^6-264 x^2 t^5+64 x^3 t^4+32 t^4-32 x t^3+9 x^2 t^2)T^2
  \\&\quad+(192 x^3 t^7+128 t^7-96 x^4 t^6-192 x t^6+128 x^2 t^5-32 x^3 t^4)T^3
  \\&\quad+(48 x^5 t^8+192 x^2 t^8-192 x^3 t^7+56 x^4 t^6)T^4
  \\&\quad+(96 x^4 t^9-48 x^5 t^8)T^5+16 x^6 t^{10}T^6.  
\end{alignat*}

Running \textsf{Maple12} on a modern laptop\footnote{MacBook Pro; Intel Core 2 Duo Processor, @2.4 GHz; 4Mb cache, 2Gb RAM.}, the whole guessing computation requires about 80Mb of memory and takes less than 20 seconds.  Once the candidate polynomial $P$ is guessed, one could proceed to its empirical certification; this can be done in various ways, as explained in~\cite{bostan09}.               
We do not need to do this here, since we are going to \emph{prove\/} in \S\ref{ssec:proof} that $F(t;x,0)$ is a root of $P$.

\medskip 
One may wonder where the precision 80 used in the previous computations
comes from. Here, this precision was humanly guessed, being 
chosen as a reasonable threshold. However, a straightforward doubling technique (not explained here in detail) would allow to \emph{automatically\/} tune it, by running several times the whole guessing procedure with increasing precision until the same polynomial is output two consecutive times.

\subsection{Proving}   \label{ssec:proof}
In this section, we detail the two steps~\eqref{step:2a} and~ 
\eqref{step:2b} used in the proof of Theorem~\ref{thm:k}.

\subsubsection{Existence and Uniqueness}

Since $P(1,0,x)=0$ and $\frac{\partial P}{\partial T} (1,0,x)=-x$, the  
implicit function theorem implies that $P$ admits a unique root  
$F_{\cand}(t;x,0)$ in $\set Q((x))[[t]]$. It follows that
$P$ has \emph{at most one\/} root in $\set Q[[x,t]]$ and that
this root, if it exists, belongs to $\set Q[x,x^{-1}][[t]]$.

Proving \emph{the existence\/} of a root of $P$ in $\set Q[[x,t]]$ is  
less straightforward: this time,
the equalities $P(1,0,0)=0$ and $\frac{\partial P}{\partial T}  
(1,0,0)=0$
prevent us from directly invoking the implicit function theorem.

\medskip We are thus faced to a clumsy technical complication, since  
what we really
need to prove is that the root $F_{\cand}(t;x,0)$ actually belongs to
$\set Q[[x,t]]$: otherwise, the substitution of $U = F_{\cand}(t;x,0)$  
in Equation~\eqref{eq:Kred}, used in Step~\eqref{step:2b} of the proof  
of Thm.~\ref{thm:k}, would not be legitimate.

\medskip
To circumvent this complication, we exploit the fact that, when seen  
in $\set Q(x)[T,t]$,
the polynomial $P(T,t,x)$ defines a curve of genus zero over $\set Q(x) 
$, which can thus be rationally parameterized. Precisely, using  
\textsf{Maple}'s \textsf{algcurves} package, the rational functions  
$R_1(U,x)$ and $R_2(U,x)$ defined by:
\begin{alignat*}1
  R_1(U,x)&=
\frac{U (1+U) (1+2 U+U^2+U^2 x)^2}{h(U,x)}, \\
R_2(U,x)&=\frac{(U^4x^2+2U^2(U+1)^2x+1+4U+6U^2+2U^3-U^4) h(U,x)}{
(1+U)^2 (1+2U+U^2+U^2x)^4
},
\end{alignat*}
with
\[h(U,x) = U^6x^3+3U^4(U+1)^2x^2+3U^2(U+1)^4x+1+6U 
+15U^2+24U^3+27U^4+18U^5+5U^6,\]
are found to share the following properties:
\begin{enumerate}
\item[$\bullet$] $P(R_2(U,x),R_1(U,x),x) = 0$;
\item[$\bullet$] there exists a (unique) power series
\[U_0(t,x) = t+t^2+(x+1)t^3+(2x+5)t^4+(2x^2+3x+9)t^5+\ldots\] in $\set  
Q[[x,t]]$ such that $R_1(U_0,x) = t$ and $U_0(0,x)=0$.
\end{enumerate}

While the first property is easily checked by direct calculation, the second  
one is a consequence of the implicit function theorem,
since $Q(U,t,x)=R_1(U,x)-t$ satisfies $Q(0,0,0)=0$ and $\frac{\partial  
Q}{\partial U}(0,0,0)=1$.

The existence proof of a power series solution of $P$ is then  
completed using the following argument: $R_2$ having no pole at $U=0$,  
and the valuation of $U_0$ with respect to $t$ being positive, the  
composed power series $R_2(U_0(t,x),x)$ is well defined in $\set  
Q[[x,t]]$ and it satisfies
\[P(R_2(U_0,x),t,x) = P(R_2(U_0,x),R_1(U_0,x),x) = 0.\]
Therefore, $F_{\cand}(t;x,0) = R_2(U_0(t,x),x)$ is the unique power  
series solution in $\set Q[[x,t]]$ of $P$.


\subsubsection{Compatibility with the reduced kernel equation}

We need to show that $F_{\cand}(t;x,0)$ so defined satisfies equation  
\eqref{eq:Kred}.
This can be done in various ways by resorting to closure properties  
for algebraic power
series. These closure properties are performed by means of resultant  
computations, based on Lemma~\ref{lem:4} below. 

\medskip 
One possibility is to first prove that the power series $S(t,x) \in \set Q[x,x^{-1}][[t]]$ defined by
\[
S(t,x) = \frac{Y(t,x)}{t} - \frac{Y(t,x)}{x} F_{\cand}(t;Y(t,x),0)
\]
is a root of the polynomial $P(T,t,x)$, and then to use the fact that $P$ 
has only one root in $\set Q[x,x^{-1}][[t]]$, namely 
$F_{\cand}(t;x,0)$. This will imply that $S(t,x)$ and $F_{\cand}(t;x,0)$ coincide, and thus that $F_{\cand}(t;x,0)$ satisfies equation \eqref{eq:Kred}, as desired.

\medskip The main point of this approach is that, since the power series $Y(t,x)$ and $F_{\cand}(t;x,0)$ are both algebraic, finding a polynomial which annihilates the series $S(t,x)$ can be done in an \emph{exact\/} manner, without having to appeal to guessing routines. Moreover, the minimal polynomial of $S(t,x)$ can be determined by factoring an annihilating polynomial obtained through a resultant computation, and, if necessary, by matching the irreducible factors against the initial terms of the series~$S(t,x)$.

More precisely, one can use the following classical facts, that we recall for completeness, see e.g.~\cite{loos83} for a proof.                        

\begin{lem}\label{lem:4}
Let $\set K$ be a field and let $P,Q\in\set K[T,t,x]$ be annihilating polynomials of two algebraic power series $A,B$ in $\set K[x,x^{-1}][[t]]$. Then
\begin{enumerate}
\item\label{scalmul} $p A$ is algebraic for every $p\in\set K(t,x)$, and it is a root of $p^{\deg_T P} P(T/p,t,x)$.
\item\label{sum} $A\pm B$ is algebraic, and it is a root of $\res_z(P(z,t,x),Q(\pm(T-z),t,x))$.
\item\label{mul} $AB$ is algebraic, and it is a root of $\res_z(P(z,t,x),z^{\deg_T Q}Q(T/z,t,x))$. 
\item\label{algsubst} If $\ord_x B>0$, then $A(t,B(t,x))$ is algebraic, and
   it is a root of $\res_z(P(T,t,z),Q(z,t,x))$.
 \end{enumerate}         
\end{lem}
                  
Since $z/t - z/x F_{\cand}(t;z,0)$ is
a root of (the numerator of)
$P(x/z(z/t-T),t,z)$ and since
$Y(t,x)$ is 
a root of 
$(x + T + x^2 T^2) t - xT$, 
Lemma~\ref{lem:4} suggests continuing our \textsf{Maple} session 
by constructing 
a polynomial in $\set Q [T,t,x]$ which has $S(t,x)$ as a root, in the following way:
\begin{verbatim}      
> ker := (T,t,x) -> (x+T+x^2*T^2)*t-x*T:  
> pol := unapply(P,T,t,x):
> res := resultant(numer(pol(x/z*(z/t-T),t,z)), ker(z,t,x), z):
> factor(primpart(res,T)); 
\end{verbatim}
The output of the last line is $P(T,t,x)^2$, 
which proves that $S(t,x)$ is a root of $P(T,t,x)$. 

\subsection{Consequences}  
Setting $x$ to $0$ in $P$ leads to the conclusion that 
the generating series  $F(t;0,0)$ of Kreweras excursions is 
a root of 
the polynomial 
$64t^6T^3+16t^3T^2+T-72t^3T+54t^3-1$.  
An argument similar to that used in the proof of Corollary~\ref{cor:2}
then implies that the coefficients $a_n$ of $F(t;0,0)$ satisfy the 
linear recursion
\[
(n + 6) (2 n + 9) a_{n + 3} - 54 (n + 2) (n + 1) a_n = 0, \qquad a_0 = 1, \, a_1 = 0, \, a_2 = 0, 
\]      
which in turn provides an alternative proof of the classical fact~\cite{kreweras65,gessel86,bousquet05} that the series $F(t;0,0)$ is both algebraic and hypergeometric, and it has the following closed form
\[F(t;0,0)=
{_3F_2}\biggl(\begin{matrix}
  {1/3}\kern.707em {2/3}\kern.707em{1}\\ 
  {3/2}\kern1em{2}
\end{matrix}\,\bigg|\,27\,t^3\biggr)
=\sum_{n=0}^\infty\frac{4^n \binom{3n}{n} }{(n+1)(2n+1)} t^{3n}.
\]


\section{Gessel walks}\label{sec:gessel}\label{sec:proof}

For establishing the proof of Theorem~\ref{thm:g}, we apply essentially the same reasoning 
that was applied in the previous section for proving Theorem~\ref{thm:k}. The main difference
is that the intermediate expressions get very big, so that they can only be handled by special
purpose software (see the data provided on our website~\cite{bostan08a}). There are also some 
additional complications which require to vary the arguments slightly. In this section,
we point out these complications, describe how to circumvent them, and we document our
computations.

The numbers $g(n,i,j)$ of Gessel walks of length~$n$ ending at $(i,j)\in\set Z^2$ satisfy
the recurrence equation
\[
 g(n+1;i,j) = g(n;i-1,j-1)+g(n;i+1,j+1)+g(n;i-1,j)+g(n;i+1,j) 
\]
for $n,i,j\geq 0$. Together with appropriate boundary conditions, this equation implies
that the generating function
\[
 G(t;x,y)=\sum_{n=0}^\infty\Bigl(\sum_{i,j=0}^\infty g(n;i,j)x^iy^j\Bigr)t^n,
\]
which we seek to prove algebraic, satisfies the equation
\begin{align}\label{eq:Kgessel}\tag{$\mathrm{K^G}$}
   ((1+y+x^2y+x^2y^2)t - xy)
    G(t;x,y)
     = (1 + y) t \,G(t; 0, y) + t\,G(t; x, 0) - t\,G(t; 0, 0) - xy.
\end{align}
This is the starting point for the kernel method. 

In this case, because of lack of symmetry with respect to $x$ and~$y$, there are two different
ways to put the left hand side to zero, using the  
two substitutions
\begin{alignat*}1
 y\to Y(t,x):={}&-\bigl(t x^2-x+t+\sqrt{(t x^2-x+t)^2-4 t^2 x^2}\bigr)\big / (2 t x^2) \\
             ={}&\tfrac1x t+
	      \tfrac{x^2+1}{x^2}t^2+
	      \tfrac{x^4+3x^2+1}{x^3}t^3+
	      \tfrac{x^6+6 x^4+6 x^2+1}{x^4}t^4+ \cdots \\
\text{and}\quad
 x\to X(t,y):={}&\bigl(y-\sqrt{y(y-4 t^2 (y+1)^2)}\bigr)\big/(2 t y (y+1)) \\
             ={}&\tfrac{y+1}{y}t+
	      \tfrac{(y+1)^3}{y^2}t^3+
	      \tfrac{2(y+1)^5}{y^3}t^5+
              \tfrac{5(y+1)^7}{y^4}t^7+\cdots
\end{alignat*}
They yield the equations
\begin{alignat*}1\label{eq:Kredgessel}\tag{$\mathrm{K_{red}^G}$}
 \begin{split} 
	   G(t;x,0)&= x Y(t,x)/t + G(t;0,0) - (1+Y(t,x))G(t;0,Y(t,x)),\\
   (1+y)G(t;0,y)&= X(t,y) y/t + G(t;0,0) - G(t;X(t,y),0),         
\end{split}    
\end{alignat*}
respectively.
Note that the first equation is free of $y$ while the second is 
free of~$x$. If we rename $y$ to $x$ in the second equation, then
all quantities belong to $\set Q[x,x^{-1}][[t]]$.
Note also that we can write $G(t;x,0)=G(t;0,0)+x U(t,x)$
and $G(t;0,x)=G(t;0,0)+ x V(t,x)$ for certain power series
$U,V\in\set Q[[x,t]]$.
In terms of $U$ and~$V$, the two equations above are then equivalent to          
\begin{alignat*}1 \label{eq:Kredgessel2}\tag{$\mathrm{K_{red}^{G,2}}$}  
 \begin{split}
       x U(t,x)&= x Y(t,x)/t - (1+Y(t,x))G(t;0,0) - Y(t,x)(1+Y(t,x)) V(t,Y(t,x)),\\
  (1+x)x V(t,x)&= X(t,x) x/t -(1+x)G(t;0,0) - X(t,x)U(t,X(t,x)).
 \end{split}
\end{alignat*}
The two equations \eqref{eq:Kredgessel2} correspond to the equation \eqref{eq:Kred} 
in Section~\ref{sec:kreweras}.
The situation here is more complicated in two respects. First, we have two equations 
and two unknown power series $U$~and~$V$ rather than a single equation with a single 
unknown power series $F(t;x,0)$; this difference originates from the lack of symmetry
of $G(t;x,y)$ with respect to $x$ and~$y$, 
which itself comes from the asymmetry of the Gessel step set with respect to the main diagonal of $\set N^2$.
Second, the two equations for $U$ and $V$
still contain $G(t;0,0)$ while there is no term $F(t;0,0)$ present in \eqref{eq:Kred};
this difference originates from the fact that Gessel's step set contains the admissible
step~$\swarrow$, 
as opposed to Kreweras's step set. 
The 
occurrence
of $G(t;0,0)$ in the equations~\eqref{eq:Kredgessel2} is not really problematic, as we know 
this power series explicitly thanks to Theorem~\ref{thm:0}. As for the other difference,
we need the following variation of Lemma~\ref{lem:1}.

\begin{lem}\label{lem:2}
 Let $A_1,A_2,B_1,B_2,Y_1,Y_2\in\set Q[x,x^{-1}][[t]]$ be 
 such that $\ord_t B_1>0$,
 $\ord_t B_2>0$, $\ord_t Y_1>0$ and $\ord_t Y_2>0$. 
 Then there exists at most one pair $(U_1,U_2)\in\set Q[[x,t]]^2$ with 
 \begin{alignat*}1
   U_1(t,x) &= A_1(t, x) + B_1(t, x) \cdot U_2(t, Y_1(t, x)),\\
   U_2(t,x) &= A_2(t, x) + B_2(t, x) \cdot U_1(t, Y_2(t, x)).
 \end{alignat*}
\end{lem}              
\begin{proof}
By linearity, it suffices to show that the only solution $(U_1,U_2)$ in $\set Q[[x,t]] \times \set Q[[x,t]]$ of the homogeneous system  \begin{alignat*}1
   U_1(t,x) &=  B_1(t, x) \cdot U_2(t, Y_1(t, x)),\\
   U_2(t,x) &=  B_2(t, x) \cdot U_1(t, Y_2(t, x))
 \end{alignat*} is the trivial solution $(U_1,U_2)=(0,0)$. This is a direct consequence of the fact that if both $U_1$ and $U_2$ were non-zero, then the valuation of $B_1(t, x) \cdot U_2(t, Y_1(t, x))$ would be strictly greater than the valuation of $U_2(t,x)$, and the valuation of $B_2(t, x) \cdot U_1(t, Y_2(t, x))$ would be strictly greater than the valuation of $U_1(t,x)$, thus $\ord_t(U_1)>\ord_t(U_2)>\ord_t(U_1)$, a contradiction. Therefore, one of $U_1,U_2$ is zero, and the system then implies that both are zero.     
\end{proof}

By a slightly more careful analysis, the lemma could be refined further such as to show
that there is only one triple of power series $(U,V,G)$ with $U,V\in\set Q[[x,t]]$ and 
$G\in\set Q[[t]]$ (free of~$x$) which satisfies~\eqref{eq:Kredgessel2} with $G(t;0,0)$ replaced by $G$. In this version, 
the proof could be completed without reference to the independent proof of Thm.~\ref{thm:0}.

Either way, we can in principle proceed from this point as in Section~\ref{sec:kreweras}.
Out of convenience, we choose to regard $G(t;0,0)$ as known.
Again, we divide the remaining task in two steps:
\begin{enumerate}
\item\label{gstep:1} \emph{Guess} defining algebraic equations for $U(t,x)$ and $V(t,x)$,
by inspecting the initial terms of $G(t;x,0)$, resp. of $G(t;0,x)$.
\item\label{gstep:2} \emph{Prove} that 
\begin{enumerate}   
           \item\label{gstep:2a}  each of the guessed equations has a unique solution in $\set Q[[x,t]]$, denoted $U_{\cand}(t;x,0)$, resp. $V_{\cand}(t;x,0)$;
 	   \item\label{gstep:2b} the power series $U_{\cand}$ and $V_{\cand}$ 
indeed satisfy the two equations in~\eqref{eq:Kredgessel2}.
\end{enumerate}
\end{enumerate}

Once this has been accomplished, Lemma~\ref{lem:2} implies that the candidate series are
actually equal to $U$ and $V$, respectively, and so these series 
as well as $G(t;x,0)$ and $G(t;0,y)$ are in 
particular algebraic. Then equation \eqref{eq:Kgessel} implies that $G(t;x,y)$ is 
algebraic, too. This then completes the proof of Thm.~\ref{thm:g}.

\subsection{Guessing}\label{sec:3.1}  

In the beginning, we had no reason to suspect that $G(t;x,y)$ is algebraic, since even the specialization $G(t;0,0)$ was generally thought to be transcendental.

Motivated by the case $x=y=0$
(i.e., by Thm.~\ref{thm:0}, which was merely a conjecture by that time), we wanted to find out whether $G(t;x,y)$ has chances to be D-finite with respect to~$t$, and searched for linear
differential equations with polynomial coefficients potentially satisfied by its sections $G(t;x,0)$ and
$G(t;0,y)$. With such equations at hand, we could have, in principle, proven the 
D-finiteness of $G(t;x,y)$,
by very much the same reasoning that we apply here for proving that $G(t;x,y)$ is algebraic.
       

We realized quickly that the differential equations for $G(t;x,0)$ and $G(t;0,y)$, if they exist, are
too big to be caught by the guessers implemented in 
packages like 
\textsf{Maple}'s \textsf{gfun} or 
\textsf{Mathematica}'s \textsf{GeneratingFunctions}.     
In order to gain 
efficiency, we switched to \textsf{Magma}, 
which provides efficient implementations of low-level algorithms,
and we opted for applying 
a modular approach:
we 
set $x$ and $y$ to special values $x_0,y_0=1,2,3,\dots$, 
and in addition, we kept numerical coefficients reduced modulo several 
fixed primes~$p$ to avoid the emerging of large rational numbers. 
Modulo a prime~$p$, and starting from the first 1000 terms of the 
series $G(t;x,0)$ and
$G(t;0,y)$, we used a very efficient automated guessing scheme, relying on the Beckermann-Labahn (FFT-based) super-fast algorithm for computing Hermite-Pad\'e approximants~\cite{beckermann94}. Eventually, we made the following observations:

  \begin{itemize}
  \item For any choice of $p$ and~$x_0$, there are several differential operators 
in $\set Z_p[t]\langle D_t\rangle$, 
of order~14 and with coefficients of degree at most~43, 
which seem to annihilate $G(t;x_0,0)$ in $\set Z_p[[t]]$.     
  \item For any choice of $p$ and~$y_0$, there are several differential operators              
in $\set Z_p[t]\langle D_t\rangle$, 
of order~15 and with coefficients of degree at most~34,  
which seem to annihilate $G(t;0,y_0)$ in $\set Z_p[[t]]$.                                                       
  \end{itemize}               

                  

\noindent (Here, and hereafter, $D_t$ stands for the usual derivation operator $\frac{d}{dt}$, and for any 
ring~$R$, we denote by $R[t]\langle D_t\rangle$  the Weyl algebra of differential operators with polynomial coefficients in $t$ over $R$.)       
                       

\medskip

The next idea was to apply an interpolation mechanism in order to 
reconstruct, starting from guesses for various choices of $x_0$ and $y_0$, and modulo various primes $p$, two candidate operators: one in $\set Q[x,t]\langle D_t\rangle$ that would annihilate $G(t;x,0)$ in $\set Q[x][[t]]$, and the other one in $\set Q[y,t]\langle D_t\rangle$, that would annihilate $G(t;0,y)$ in $\set Q[y][[t]]$. The ingredients needed to put such an interpolation scheme into practice are \emph{rational function interpolation\/} and \emph{rational number reconstruction\/}.
Both are standard techniques in computer algebra, for details on fast algorithms we refer to \cite{vzgathen99}.

\medskip

To our surprise, when applied to the order 14 and 15 differential operators mentioned above, we found this reconstruction scheme to require an unreasonably large number
of evaluation points $x_0,y_0=1,2,3,\dots$, suggesting an unreasonably high degree of the
operators with respect to $x$ or~$y$, respectively. We aborted the computation when 
the expected degree exceeded~{}1500 (!). 
At this point, we had the impression that the section series $G(t;x,0)$ and $G(t;0,y)$ might not be D-finite.

\medskip
Our next 
attempt was to find candidate operators of smaller total 
size by trading order against degree. We went back to the series $G(t;x_0,0)$ modulo $p$, 
and 
tried to determine the \emph{least order\/} operator $\mathcal{L}_{x_0,0}^{(p)} \in \set Z_p[t]\langle D_t\rangle$ annihilating it.  This was done by taking several candidate operators in $\set Z_p[t]\langle D_t\rangle$ of order~14 as above, and by computing their greatest common right divisor (gcrd) in the rational Weyl algebra $\set
Z_p(t)\langle D_t\rangle$.  (Despite the non-commutativity of $\set Z_p(t)\langle D_t \rangle$, this can be
done by a variant of the Euclidean algorithm~\cite{ore33,bronstein96}, see~\cite{grigoriev90} for a more efficient grcd algorithm.) 
We applied the same strategy to find the \emph{least order\/} operator $\mathcal{L}_{0,y_0}^{(p)} \in \set Z_p[t]\langle D_t\rangle$ annihilating 
the series $G(t;0,y_0)$ modulo $p$.

Doing so for
several evaluation points $x_0,y_0=1,2,3,\dots$ and several choices of~$p$, it was
finally possible, by using the interpolation scheme described above, to reconstruct from  
the various modular candidate operators $\mathcal{L}_{x_0,0}^{(p)}$ and $\mathcal{L}_{0,y_0}^{(p)}$, 
two candidates
$\mathcal{L}_{x,0}  \in \set Q[x,t]\langle D_t\rangle$ and $\mathcal{L}_{0,y} \in \set Q[y,t]\langle D_t\rangle$ with
\emph{reasonable\/} degrees in $x$ and~$y$.   

The operators $\mathcal{L}_{x,0}$ and $\mathcal{L}_{0,y}$ are posted on our
website~\cite{bostan08a}. 
The operator $\mathcal{L}_{x,0}$ 
has order~11, degree~96
with respect to~$t$ and degree only 78\,(!) in~$x$. Its longest integer
coefficients has only 61 decimal digits.  The operator $\mathcal{L}_{0,y}$ 
is even
nicer. Its order is~11, its degree is~68 with respect to~$t$ and just 28\,(!)
with respect to~$y$. Its longest integer coefficient has 51 decimal digits.

\medskip
The whole procedure for guessing $\mathcal{L}_{x,0}$ and $\mathcal{L}_{0,y}$ 
took less than 2 CPU hours on a modern computer running \textsf{Magma v 2.13 (12)}.
To achieve this speed, we greatly benefited, 
on the one side, from the fast \textsf{Magma}'s buit-in polynomial, integer and modular arithmetic, and on the other side, from our own efficient implementations of several algorithms (e.g., for Hermite-Pad\'e approximation, for rational function interpolation and for grcds).  

For instance, to compute $\mathcal{L}_{x,0}$ we used 21 primes $p$ of 28 bits each and 158 distinct integer values of~$x_0$. 
Modulo each prime $p$, 150 CPU seconds were enough to compute:
$158$ bundles of four differential operators in $\set Z_p[t]\langle D_t\rangle$ with order $14$ and with coefficients of degree at most $43$ (by Hermite-Pad\'e approximation), 158 operators in $\set Z_p[t]\langle D_t\rangle$ of order $11$ and with coefficients of degree at most $96$  (by right grcd)  and $12 \times 97 = 1164$ rational functions  in $\set Z_p(x)$  with numerators and denominators of degree at most 78 (by rational interpolation). At this point, the operator contained 
$97 \times 79 \times 12 = 91956$ terms of the form $c_{i,j,k} \,t^i x^j D_t^k$  ($i \leq 96, j \leq 78, k \leq 11$), where each $c_{i,j,k} \in \set Q$ was known modulo the 21 primes. The constants $c_{i,j,k}$ were recovered by performing 91956 rational number reconstructions. The whole computation of $\mathcal{L}_{x,0}$ took 55 minutes.
Guessing $\mathcal{L}_{0,y}$ using the same method was even a little bit faster.

\medskip The exceptionally small sizes of $\mathcal{L}_{x,0}$ and $\mathcal{L}_{0,y}$ (in comparison to the intermediate expressions)  speak
very much in favor of their correctness. 
Also, 
the fact that 
the operators $\mathcal{L}_{x,0}$ and $\mathcal{L}_{0,y}$ verify the following equalities in $\set Q[[t]]$:  
\[\mathcal{L}_{x,0} (G(t;x,0)) \bmod t^{1000} =  0 \quad \text{and} \quad      \mathcal{L}_{0,y}(G(t;0,y)) \bmod t^{1000} = 0,\]
provides more empirical evidence that $\mathcal{L}_{x,0}$ and $\mathcal{L}_{0,y}$ 
are indeed annihilating operators for $G(t;x,0)$ and $G(t;0,y)$, respectively.

There are a number of additional
tests which can be performed to experimentally sustain the evidence that a guessed differential operator
is correct (see our paper~\cite{bostan09} for a collection of such tests), and our operators $\mathcal{L}_{x,0}$ and $\mathcal{L}_{0,y}$ successfully pass 
all these tests. 

\medskip         

One of the tests consists of checking whether the operators $\mathcal{L}_{x,0}$ and $\mathcal{L}_{0,y}$ possess an arithmetic property which is expected from the minimal order operator annihilating 
a generating function like $G(t;x,0)$ and $G(t;0,y)$, see~\cite{bostan09}. This property, called \emph{global nilpotency\/} \cite{dwork90}, can be stated as follows: 
for almost any prime number~$p$, the order~11 operators $\mathcal{L}_{x,0}$, resp. $\mathcal{L}_{0,y}$, should right-divide the pure power $D_t^{11\cdot p}$ in $\set Z_p(x,t)\langle D_t\rangle$, resp. in $\set Z_p(y,t)\langle D_t\rangle$. We checked that this property indeed holds for all primes $p < 100$. We actually found out that  
the operators $\mathcal{L}_{x,0}$ and $\mathcal{L}_{0,y}$ have a stronger property: they even right-divide $D_t^{p}$; in other terms, they have zero $p$-curvature for all the tested primes~$p$. This was the key observation 
which led us to suspect that $G(t;x,y)$ is algebraic, for according to a famous conjecture
of Grothendieck~\cite{put01}, an operator has zero $p$-curvature if and only if it admits a basis of algebraic solutions.
The conjecture is still open (even for second order operators), but
it is generally believed to be true. In either case, something interesting is going on:
either $G(t;x,y)$ is algebraic or we have found operators which very much look like 
counterexamples to Grothendieck's conjecture.

\medskip We next searched for potential polynomial equations satisfied by
the power series $U,V\in\set Q[[x,t]]$ defined by 
$G(t;x,0)=G(t;0,0)+x U(t,x)$
and $G(t;0,x)=G(t;0,0)+ x V(t,x)$.
We did not find any using only 1000 terms of those series, but
we found some starting from 1200 terms. Using again guessing techniques based on fast modular Hermite-Pad\'e approximation, combined with an interpolation scheme, 
we discovered two polynomials $P_1(T,t,x) \in \set Q[T,t,x]$ and $P_2(T,t,y)\in \set Q[T,t,y]$
which 
satisfy
\[
 P_1(U(t,x),t,x)=0 \bmod {t^{1200}}\quad\text{and}\quad P_2(V(t,y),t,y) = 0 \bmod {t^{1200}}.
\]
These polynomials are posted on our website~\cite{bostan08a}. The polynomial $P_1$ has
degrees 24,~44, and~32 with respect to $T$, $t$, and~$x$, respectively, and involves 
integers with no more than 21 decimal digits. The polynomial $P_2$ has degrees 24,~46, and~56
with respect to $T$, $t$, and~$y$, respectively, and involves integers with no more 
than 27 decimal digits. Spelled out explicitly in this article, they would both together 
fill about thirty pages; they are however much smaller than the differential operators $\mathcal{L}_{x,0}$ and $\mathcal{L}_{0,y}$, for which five hundred pages would not be enough!
Just like 
$\mathcal{L}_{x,0}$ and $\mathcal{L}_{0,y}$, the polynomials $P_1$ and $P_2$ pass
a number of heuristic tests which let them appear plausible. 

\medskip 
We are now going to prove that the guessed polynomials $P_1$ and $P_2$ are indeed valid.

\subsection{Proving}

 Let $P_1\in\set Q[T,t,x]$ and $P_2\in\set Q[T,t,y]$ be the two polynomials
 posted on the website to this article~\cite{bostan08a}. 
 We show (i) that these polynomials admit unique power series solutions
 $U_{\cand}(t,x)$ and $V_{\cand}(t,x)$, respectively, and (ii) that these
 power series satisfy the reduced kernel equations~\eqref{eq:Kredgessel2}.

 \subsubsection{Existence and Uniqueness}

 As in the case of Kreweras's walks, the implicit function theorem does not
 apply to these polynomials, but unlike in the Kreweras case, an existence 
 proof using a suitable rational parameterization is not possible either, 
 because the polynomials at hand define curves of positive genus, and 
 therefore a rational parameterization does not exist. 

 In order to obtain a proof in this situation, we proceeded along the following 
 lines:
 \begin{itemize}
 \item First we used Theorem~3.6 of McDonald~\cite{mcdonald95} to
   obtain the existence of a series solution
   \[
     \sum_{p,q\in\set Q} c_{p,q}t^px^q
   \]
   with $c_{p,q}=0$ for all $(p,q)$ outside a certain halfplane~$H\subseteq\set Q^2$.
 \item Next, we computed a system of bivariate recurrence equations
   with polynomial coefficients that the coefficients $c_{p,q}$ must
   necessarily satisfy. This can be done in principle by software packages
   such as Chyzak's \texttt{mgfun}~\cite{chyzak98} or Koutschan's 
   \texttt{HolonomicFunctions.m}~\cite{koutschan09}.
   However, for reasons of efficiency we used our own implementation of
   the respective algorithms.
 \item The form of the recurrences together with the shape of the halfplane~$H$ imply 
   that the coefficients $c_{p,q}$ of any solution can be nonzero only in
   a finite union of cones $v+\set N u+\set N w$ with vertices $v\in\set Q^2$ 
   and basis vectors $u,w\in\set Q^2$ that can be computed explicitly. 
   If $c_{p,q}\neq0$ for some index $(p,q)$ in such a cone, then also the 
   coefficient at the cone's vertex must be nonzero.
 \item Applying McDonald's generalization of Puiseux's algorithm, we 
   determined the first coefficients of series solutions 
   to an accuracy that all further coefficients belong to some translate
   of $H$ which contains no vertices. 
 \item As one of these partial solutions contained no terms with fractional
   powers, it was possible to conclude that the entire series contains
   no terms with fractional exponents. Reference to $u$ and~$w$ implied that 
   this partial solution could also not contain any terms with negative 
   integral exponents,
   so the only remaining possibility was that the solution is in fact a 
   power series. 
 \end{itemize}

 A full description of the argument requires a somewhat lengthy
 discussion of a number of technical details, which we prefer to avoid here. 
 A supplement to this article is provided on our website~\cite{bostan08a}
 in which we carry out existence proofs in full detail that
 both $P_1$ and $P_2$ admit some power series solutions $U_{\cand}$ and $V_{\cand}$,
 respectively. 

 \subsubsection{Compatibility with the reduced kernel equation}

 It remains to show that these
 solutions $U_{\cand}$ and $V_{\cand}$ satisfy the system~\eqref{eq:Kredgessel2}.
 Because of $X(t,Y(t,x))=x$, the substitution $x\to Y(t,x)$ transforms the
 second equation of that system to the first. 
 Therefore, it suffices to prove the second equation:
 \begin{equation} \label{eq:second}
  (1+x)x V_{\cand}(t;x,0)= X(t,x) x/t -(1+x)G(t;0,0) - X(t,x)U_{\cand}(t;X(t,x),0).
 \end{equation}

Letting $G_1(t,x) = G(t;0,0) + xU_{\cand}(t;x,0)$ and $G_2(t,x) = G(t;0,0) + xV_{\cand}(t;x,0)$, the last equation is equivalent to 
\begin{equation} \label{eq:secondPrime}                                     
 (1+x) G_2(t,x) - G(t;0,0) = x X(t,x)/t - G_1(t,X(t,x)).
\end{equation}


 By Corollary~\ref{cor:2} and Lemma~\ref{lem:4}, the power series 
 \[
   (1+x) G_2(t,x) - G(t;0,0)
   \quad\text{and}\quad
   x X(t,x)/t - G_1(t,X(t,x))
 \]
 are algebraic and we can compute their minimal polynomials---at least in theory.
 Now the polynomials $P_1$ and $P_2$ are so big that the required resultant 
 computations cannot be carried out by \textsf{Maple} or \textsf{Mathematica}. 

 There are efficient special purpose algorithms available for the particular kind of resultants 
 at hand~\cite{bostan06} and our \textsf{Magma} implementation of these algorithms 
is able to perform the necessary computations. It turns out that 
 the minimal polynomials for both power series are identical. 
 It is provided electronically on the website to this article.
 After determining a suitable number of initial terms of both series and observing that
 they match, it can be concluded that Equations~\eqref{eq:secondPrime} and~\eqref{eq:second} hold.
 This completes the proof of Theorem~\ref{thm:g}.

\subsection{Consequences}

The fact that $G(t;x,y)$ is algebraic has consequences which are of combinatorial interest.
We list some. 


\begin{cor}
 The following series are algebraic:
 \begin{itemize}
 \item $G(t;1,1)$ -- the generating function of Gessel walks with arbitrary endpoint.
 \item $G(t;1,0)$ and $G(t;0,1)$ -- the generating functions of Gessel walks ending
   somewhere on the $x$-axis or the $y$-axis, respectively.
 \end{itemize}
\end{cor}

Using the built-in equation solver of Mathematica, we found that
all these series, as well as the series $G(t;0,0)$, can be expressed in terms of 
nested radicals, for example
\begin{alignat*}1
 G(t;1,1)&=\frac1{6t}\left(-3+\sqrt{3}\sqrt{U(t)+\sqrt{\frac{16 t (2 t+3)+2}{(1-4 t)^2 U(t)}-U(t)^2+3}}\right)\\
 \quad\text{where}\quad
 U(t)&=\sqrt{1+4 \sqrt[3]{t(4 t+1)^2/(4 t-1)^4}}.
\end{alignat*}
The radical representations of the other series are much more involved than this one.
They are available electronically at the website to this article~\cite{bostan08a}.
Also their minimal polynomials can be found there. 


 \begin{cor}
  For every point $(i,j)$, the series 
  $G_{i,j}(t):=\sum_{n=0}^\infty g(n;i,j)t^n$ is algebraic. 
 \end{cor}
 \begin{proof}
  We have
  \[
   G_{i,j}(t) = \frac1{i!j!}\Bigl(\frac{d^i}{dx^i}\frac{d^j}{dy^j} G(t;x,y)\Bigr)\Bigr|_{x=y=0}
  \]
  and the property of being algebraic is preserved under differentiation and evaluation.
 \end{proof}

 The previous corollary implies in particular that the conjecture of Petkov{\v s}ek and 
 Wilf~\cite{petkovsek08} that
 $g(n;0,j)$ (for fixed~$j$) and $g(n;1,0)$ are P-finite 
 is right, and that their 
 conjecture that $g(n;2,0)$ is not P-finite is wrong.

 For the degrees of the minimal polynomials $p_{i,j}(T,t)$ of $G_{i,j}(t)$, we observed
 empirically that
 \begin{alignat*}1
   \deg_T p_{i,j} &= \left\{\begin{array}{ll}
        4 & \text{if $i=2j+1$}\\
        8 & \text{otherwise}
      \end{array}\right.\\
  \quad\text{and } \deg_t p_{i,j} &= \left\{\begin{array}{ll}
        12j - 5i + 14 & \text{if $i\leq j$}\\
        5i + 2j +14 & \text{if $j < i < 2j+1$}\\
        6j + 9 & \text{if $i  = 2j+1$}\\
        7i - 2j + 12 & \text{if $i > 2j+1$}
      \end{array}\right.
 \end{alignat*}
 but we are not able to prove these degree formulas in general. 

 \medskip
 Note that our proof of Theorem~\ref{thm:g} does not provide us with the minimal polynomial of
 $G(t;x,y)$. This polynomial will, in fact, be much larger than the minimal polynomials
 of the series $G_{i,j}(t)$ or the series obtained form $G(t;x,y)$ by setting $x$, $y$~or~$t$
 to special values. 
 From the sizes of the minimal polynomials of $G(t;x,0)$ and $G(t;0,y)$, which we know
 explicitly, it can be deduced that the minimal polynomial $p(T,t,x,y)$ of $G(t;x,y)$
 will have degrees 72,~141, 263, and~287 with respect to
 $T$,~$t$, $x$, and~$y$, respectively, and
 thus consist of more than 750 Mio terms.

 \begin{cor}\label{cor:10}
   $G(t;x,y)$ is D-finite with respect to any of the variables~$x,y$ and $t$.
 \end{cor}

 As every algebraic power series is D-finite, this is an immediate consequence
 of Theorem~\ref{thm:g}, even if we regard $G(t;x,y)$ as a multivariate power
 series in $t$, $x$, and~$y$ rather than as a power series in $t$ only with
 $x$ and $y$ belonging to the coefficient domain. 

 D-finiteness in $t$ only amounts to the existence of a linear differential equation 
 in $d/dt$ with coefficients in $\set Q(x,y)[t]$. This can be proven independently in
 as similar way as Theorem~\ref{thm:g}.
 It suffices to discover differential operators which potentially 
 annihilate $U(t,x)$ and~$V(t,y)$, respectively, define 
 $U_{\cand}(t,x)$ and $V_{\cand}(t,y)$ as
 the unique power series annihilated by these operators and 
 matching the first terms of $U(t,x)$ and~$V(t,y)$, respectively, 
 and prove that these series satisfy the equations in \eqref{eq:Kredgessel}.

 The option of proving D-finiteness (with respect to~$t$) directly is important when other step 
 sets instead of Gessel's $\{\leftarrow,\rightarrow,\nearrow,\swarrow\}$ 
 are considered, for which the generating function $G(t;x,y)$ is D-finite in $t$ 
 but not algebraic. As shown by Mishna~
\cite{mishna09}, such step sets do exist.


 Explicit knowledge of differential operators annihilating $G(t;x,0)$
 and $G(t;0,y)$ also allows to deduce bounds on the size of the 
 differential operator annihilating~$G(t;x,y)$. According to a priori
 estimations, this operators will have order up to~22 and polynomial
 coefficients with 
 degrees up to 1968,~936, and~336 with respect to $t$, $x$, and~$y$,
 respectively, and thus consist of about $1.4\cdot10^{10}$ terms. 

 \begin{cor}\label{cor:fast}
  For fixed $i$ and~$j$, the number $g(n;i,j)$ can be computed with $O(n)$ arithmetic operations. 

  For fixed $x$ and~$y$, the coefficient $\<t^n> G(t;x,y)$ can be computed with $O(n)$ arithmetic operations. 
 \end{cor}
 \begin{proof}
  By Cor.~\ref{cor:10}, the coefficient sequence $g(n;i,j)$ is P-finite with respect to~$n$. 
  Therefore, it satisfies a uniform recurrence with respect to~$n$.
  This recurrence, together with appropriate initial values, allows the computation 
  of $g(n;i_0,j_0)$ in linear time.

  The argument for the second assertion is similar. 
 \end{proof}

\smallskip\noindent{\bf Acknowledgments.} We wish to acknowledge financial
support from the joint Inria-Microsoft Research Centre, and the Austrian Science Foundation (FWF) grants P19462-N18 and P20162-N18. 
We thank Fr\'ed\'eric Chyzak, Philippe Flajolet, Christoph Koutschan, 
Mireille Bousquet-M\'elou, Preda Mih\u{a}ilescu, Tanguy Rivoal,
Bruno Salvy, Josef Schicho and Doron Zeilberger,
for stimulating discussions during the preparation of this work. 
 \bibliographystyle{plain}
 \bibliography{gessel}

\end{document}